\documentclass[12pt]{amsart}
\usepackage{amssymb,amsfonts,latexsym}
\usepackage{graphics,verbatim}
\usepackage{graphicx}
\usepackage[usenames]{color} %\color{Red}

\newtheorem{theo}{Theorem}[section]
\newtheorem{lem}[theo]{Lemma}
\newtheorem{prop}[theo]{Proposition}
\newtheorem{coro}[theo]{Corollary}

\theoremstyle{definition}

\newtheorem{example}[theo]{Example}

\theoremstyle{remark}

\numberwithin{equation}{section}

\newcommand{\R}{{\mathbb R}}
\newcommand{\Z}{{\mathbb Z}}

\newcommand{\N}{{\mathbb N}}

\newcommand{\A}{{\mathcal A}}
\newcommand{\B}{{\mathcal B}}
\newcommand{\C}{{\mathcal C}}

\newcommand{\eproof}{\hfill$\square$}%\rule{2.2mm}{3.0mm}}

\usepackage[usenames]{color} %\color{Red}

\begin{document}

\title[Exponential spectrum]{Exponential spectra in $L^2(\mu)$}

\author[X.-G. He]{Xing-Gang He}

\address{College of Mathematics and Statistics, Central China Normal University, Wuhan 430079, China}
\email{xingganghe@yahoo.com.cn}
%\thanks{}

\author{Chun-kit Lai}
\address{Department of Mathematics , The Chinese University of Hong Kong , Hong Kong}
\email{cklai@@math.cuhk.edu.hk}
\author[K.-S. Lau]{Ka-Sing Lau}
\address{Department of Mathematics \\ The Chinese University of Hong Kong \\Shatin, Hong Kong.} \email{kslau@math.cuhk.edu.hk}
\thanks{The research is partially supported by the RGC grant of Hong Kong and the Focused
Investment Scheme of CUHK; the first author is also supported by the National Natural Science Foundation of
China 10771082 and 10871180.}

%    General info

\date{\today}
%\date{December 1, 2001 and, in revised form, June 22, 2001.}



\begin{abstract} \ Let $\mu$ be a Borel probability measure with compact support. We consider exponential type orthonormal bases, Riesz bases and frames in $L^2(\mu)$. We show that if $L^2(\mu)$ admits an exponential frame, then $\mu$ must be of pure type.  We also classify various $\mu$ that admits either kind of  exponential bases, in particular, the discrete measures and their connection with integer tiles.  By using this and convolution, we construct a class of singularly continuous  measures that has an exponential Riesz basis but no exponential orthonormal basis. It is the first of such kind of examples.

\end{abstract}

\maketitle

\section{Introduction}
\setcounter{equation}{0}

Throughout the paper we assume that  $\mu$ is a  (Borel) probability measure on ${\mathbb R}^d$ with compact support. We call a family $E(\Lambda)=\{e^{2\pi i\lambda x}: \lambda\in
\Lambda\}$ ($\Lambda$ is a countable set)   a {\it Fourier frame} of the Hilbert space $L^2(\mu)$ if there
exist $A, B >0$ such that
\begin{equation} \label {eq1.1}
A\|f\|^2 \le \sum_{\lambda\in \Lambda}|\langle f, e^{2\pi i \lambda x}\rangle|^2\le B\|f\|^2 ,  \qquad \forall \ f \in L^2(\mu).
\end{equation}
Here the inner product is defined as usual,
$$
\langle f, e^{2\pi i \lambda x}\rangle =\int_{{\mathbb R}^d}f(x)e^{-2\pi i \lambda x}d\mu(x).
$$
$E(\Lambda)$ is called an {\it (exponential) Riesz basis} if it is
both a basis and a frame of $L^2(\mu)$. Fourier frames and
exponential Riesz bases are natural generalizations of exponential orthonromal
bases in $L^2(\mu)$. They have  fundamental importance in non-harmonic
Fourier analysis  and  close connection with time-frequency analysis ([Chr], [G], [H]).  When (\ref {eq1.1}) is satisfied, $f \in L^2(\mu)$ can be expressed as $ f (x) = \sum _{\lambda \in \Lambda} c_\lambda e^{2\pi i \lambda x}$, and the expression is unique if it is a Riesz basis.

\medskip

When  $E(\Lambda)$ is an orthonormal basis (Riesz
basis, or frame) of $L^2(\mu)$, we say that $\mu$ is a {\it spectral measure (R-spectral measure, or F-spectral measure respectively)} and $\Lambda$ is called a
\textit{spectrum (R-spectrum, or F-spectrum respectively)} of $L^2(\mu)$. We will also use the term orthonormal spectrum instead of spectrum when we need to emphasis the orthonormal property.  If $E(\Lambda)$ only satisfies the upper bound condition in (\ref {eq1.1}), then it is called a \textit{Bessel set} (or \textit{Bessel sequence}); for convenience, we also call $\Lambda$  a Bessel set  of $L^2(\mu)$.

\bigskip

   One of the interesting and basic questions in  non-harmonic Fourier analysis is :

   \vspace {0.2cm}

    \textit{What kind of compactly supported probability measures in ${\mathbb R}^d$ belong to the above classes of measures? }

    \vspace {0.2cm}

  When $\mu$ is the restriction of the Lebesgue measure on $K$ with positive measure, the question whether it is a spectral measure is related to the well known Fuglede problem of translational tiles  (see [Fu], [LW], [{\L}a], [T] and the reference therein). While it is easy to show that such $\mu$ is an F-measure, it is an open question whether it is an R-spectral measure. If $K$ is an unit interval, its F-spectrum was completely classified in terms of de Brange's theory of entire functions [OS]. In another general situation,  Lai [L] proved a sharp result that if $\mu$ is absolutely continuous with respect to the Lebesgue measure, then it is an F-spectral measure if and only if its density function is essentially  bounded above and below on the support.

 \bigskip

 The problem becomes more intriguing when $\mu$ is singular. The first example of such spectral measures was given by Jorgensen and Pedersen [JP]. They showed that the Cantor measures with even contraction ratio ($\rho =  1/{2k}$) is spectral,  but the one with odd contraction ratio ($\rho = 1/{(2k+1)}$) is not.  This raises the very interesting question on the existence of an exponential Riesz basis or a Fourier frame for such measures, and more generally for the self-similar measures ([{\L}aW1,2], [DJ], [St], [HL]). In particular Dutkay {\it et al} proposed to use the \textit{Beurling dimension} as some general criteria for the existence of Fourier frame [DHSW]. They also attempted to find  a self-similar measure which admits an exponential Riesz basis or a  Fourier frame but not an exponential orthonormal basis [DHW]. However, no such examples have been found up to now.

\bigskip

In this paper,  we will carry out a  detail study of the three classes of spectra mentioned. It is known that a spectral measure must be either purely discrete or purely continuous  [{\L}aW2].  Our first theorem is a pure type law for the F-spectral measures.

\bigskip

\begin{theo} \label{th1.1} Let $\mu$ be an F-spectral measure on ${\Bbb R}^d$. Then it must be one of the three pure types:  discrete (and finite), singularly continuous or absolutely continuous.
\end{theo}

\bigskip

For the proof, the discrete case is based on the frame inequality, and the two continuous cases make use the concept of {\it lower Beurling density} of the F-spectrum.

\medskip

To complete the previous digression on  the continuous measures, we have the following conclusions for finite discrete measures.

\medskip

\begin{theo} \label{th1.2} Let $\mu = \sum_{c\in C}p_c \delta_c$ be a discrete probability measure in ${\Bbb R}^d$ with $\C$ a finite set. Then $\mu$ is an R-spectral measure.
\end{theo}

\bigskip

To determine such discrete $\mu$ to be a spectral measure, we will restrict our consideration on ${\Bbb R}^1$ and let  $\C \subset {\Bbb Z}^+$ with $0\in \C$. Then the Fourier transform of $\mu$ is
$$
\widehat \mu (x) = p_0+ p_1e^{2 \pi i c_1} x +\cdots + p_{k-1}e^{2\pi i c_{k-1}x} : =m_\mu (x),
$$
 where $P=\{p_i\}_{i=0}^{k-1}$ is a set of probability
weights. We call $m_\mu (x)$ the {\it mask polynomial} of $\mu$. Let ${\mathcal Z_\mu}= \{x \in [0,1):  m_\mu(x) =0\}$ be the zero set of $m_\mu (x)$, and $\Lambda$ is called a {\it bi-zero set} if $\Lambda -\Lambda \subset {\mathcal Z}_\mu\cup\{0\}$. Denote the
cardinality of $E$ by $\#E$. It is easy to see the following simple proposition.

\bigskip

\begin{prop} \label{th1.3} Let $\mu = \sum_{c\in \C}p_c\delta_c$ with $\C \in {\Bbb Z}^+$ and $0\in \C$.  Then $\mu$ is a spectral measure if and only if there is a bi-zero set $\Lambda$ of $m_\mu$ and $\#\Lambda = \#\C$.
In this case, all the $p_c$ are equal.
\end{prop}

\bigskip

The determination of the bi-zero set is, however, non-trivial, as the zeros of a mask polynomial is rather hard to handle.
As an implementation of the proposition, we work out explicit expressions of the set $\C$ and the bi-zero set when  $\#\C = 3, 4$. It is difficult to have such expression beyond $4$ directly. On the other hand, there are systematic studies of the zeros of the mask polynomials by factorizing the mask polynomial as cyclotomic polynomials (the minimal polynomial of the root of unity).  This has been used to study the integer tiles and their spectra (see [CM], [{\L}a], [LLR]).  We adopt this approach to a class of self-similar measures (which is continuous) in our consideration:

\medskip

Let $n>0$ and let $\A \subset {\Bbb Z}^+$ be a finite set with $0 \in {\mathcal A}$, we define a self-similar measure $\mu := \mu_{\A,n}$ by
$$
\mu(E) = \frac 1{\# A} \sum_{a \in \A} \mu (nE -a)
$$
where $E $ is a Borel subset in ${\Bbb R}$. Note that the Lebesgue measure on $[0,1]$ and the Cantor measures are such kind of measures. The following theorem is a combination of the results in [PW], [{\L}a] and [{\L}aW1]:

\bigskip

\begin{theo} \label{th1.4} \ Let $\A \subset {\Bbb Z}^+$ be a finite set with $0\in \A$. Suppose there exists $\B \subset {\Bbb Z}^+$ such that $\A \oplus \B = {\Bbb N}_n$ where ${\Bbb N}_n = \{0, \cdots , n-1\}$.  Then $\delta_{\A} = \sum_{a\in \A} \delta_a $ is a spectral measure with a spectrum in $\frac{1}{n}{\Bbb Z}$; the associated self-similar measure $ \mu_{\A,n}$ is also a spectral measure, and it has a spectrum in ${\Bbb Z}$ if gcd $\A=1$.
\end{theo}

\bigskip

Note that the $1/4$-Cantor measure $\mu_{\{0,2\},4}$ satisfies the above condition, but not the $1/3$-Cantor measure. In fact, it is an open problem whether the $1/3$-Cantor measure is an F-spectral measure. To a lesser degree we want to know the existence of a singularly continuous measure that admits an R-spectrum but is not a spectral measure. Our final goal is to search for new  R-spectral measures and to obtain such an example as corollary.

\bigskip

 To this end, we let $\eta$ be a discrete probability measure with support  $\C \subset {\Bbb Z}^+$. Let $\nu$ be another probability measure on ${\mathbb R}$ with support $\Omega\subseteq [0,1]$, and let $\mu=\eta\ast \nu$ be the convolution of $\eta$ and $\nu$. Our main result is

\bigskip

\begin{theo}\label{th1.5} Let $\mu = \eta \ast \nu$ be as the above,  and  assume that $\nu$ is an R-spectral measure with a spectrum in $ {\Bbb Z}$. Then \ $\mu$ is an R-spectral measure.

In addition, if ${\mathcal Z}_\nu\subseteq \Z$. Then $\mu$ is a spectral measure if and only if both $\eta$ and $\nu$ are spectral measures.
\end{theo}

\bigskip

We can modify the theorem slightly with the spectrum $\Gamma$ and ${\mathcal Z}_\nu$ to be some subsets of rationals (Theorems \ref{th5.2}, \ref {th5.3}), this covers some more interesting cases (e.g., the Cantor measures). Finally  by taking $\eta$ to be a {\it non-uniform} discrete measure (Proposition \ref{th1.3}) and $\nu = \mu_{\mathcal A, n}$ in Theorem \ref {th1.4}, we conclude from Theorem \ref {th1.5} that

\medskip

\begin{example} \label{th1.6} {\it There exists a singularly continuous measure which is an R-spectral measure, but not a spectral measure.}
\end{example}

\bigskip

For the organization of the paper, we prove Theorem \ref{th1.1} in Section 2 and Theorem \ref{th1.2} in Section 3. We then deal with the discrete spectral measures in Section 3;  Proposition \ref{th1.3} is proved, and explicit expressions of $\C$ (with $\#\C =3,4$) for $\mu_\C$ to be a spectral measure (Examples \ref{th3.5},  \ref{th3.6}) are sought.  In Section 4, we make a further discussion of the discrete spectral measures in connection with the class of integer tiles. In Section 5, we prove the two statements in Theorem \ref {th1.5}  in two theorems, and Example \ref{th1.6} follows as a corollary. An Appendix is included for a number-theoretic proof of Theorem \ref{th1.4}.

\bigskip
\bigskip

\section {Law of pure types}
\setcounter{equation}{0}

Recall that a  $\sigma$-finite  Borel measure $\mu$ on ${\Bbb R}^d$ can  be decomposed uniquely as \textit{discrete}, \textit{singularly continuous} and  \textit{absolutely continuous} measures, i.e., $\mu = \mu_d+\mu_s+\mu_a$.  The measure $\mu$ is said to be of \textit{pure type} if $\mu$ equals only one of the three components.

\medskip

In our proof of the pure type property of the F-spectral measures, we need to use the  \textit{lower Beurling density} of an infinite discrete set $\Lambda\subset {\Bbb R}^d$:
$$
D^{-}\Lambda := \liminf_{h\rightarrow\infty}\inf_{x\in{\Bbb R}^d}\frac{\#(\Lambda\cap Q_h(x))}{h^d},
$$
where $Q_h(x)$ is the standard cube of side length $h$ centered at $x$. Intuitively $\Lambda$ is distributed like a lattice if $D^{-}\Lambda$ is positive.  In the seminal paper [Lan], Landau gave an elegant and useful necessary condition for $\Lambda$ to be an F-spectrum on $L^2(\Omega)$: $D^{-}\Lambda\geq{\mathcal L}(\Omega)$ where ${\mathcal L}$ is the Lebesgue measure. The following proposition provides some relationships between the lower Beurling density and the types of the measures.

\medskip

\begin{prop}\label{th2.1}
Let $\mu$ be a compactly supported probability measure on ${\Bbb R}^d$ , and $\Lambda$ is an F-spectrum of $\mu$, we have

\vspace {0.15cm}
\ \ {\rm (i)} If $\mu = \sum_{c\in{\mathcal C}}p_c\delta_c$ is  discrete, then $\#\Lambda<\infty$ and $\#{\mathcal C}<\infty$;

\vspace {0.15cm}

\ {\rm (ii)} If $\mu$ is singularly continuous, then $D^{-}\Lambda =0$;

\vspace {0.15cm}

{\rm (iii)} If $\mu$ is  absolutely continuous, then $D^{-}\Lambda >0$.
\end{prop}

\medskip

\noindent{\bf Proof.} (i) \  By the definition of Fourier frame, we have for all $f\in L^2(\mu)$,
$$ \sum_{\lambda\in\Lambda}|\sum_{c\in\C}f(c)e^{2\pi i \langle\lambda, c\rangle}p_c|^2\leq B\sum_{c\in\C}|f(c)|^2p_c.
$$
 Taking $f = \chi_{c_0}$, where $p_{c_0}>0$, we have $(\#\Lambda)\cdot p_{c_0}^2\leq Bp_{c_0}$. Hence $\#\Lambda\leq B/p_{c_0} < \infty$. This implies $\#\C<\infty$ by the completeness of Fourier frame.

\bigskip

\noindent(ii) \ Suppose on the contrary that $D^{-}\Lambda\geq c>0$. We claim that  ${\Bbb Z}^d$ is a Bessel set of $L^2(\mu)$. By the definition of $D^{-}\Lambda$, we can choose a large  $h\in{\Bbb N}$  such that
$$
\inf_{x\in{\Bbb R}^d}(\#(\Lambda\cap Q_h(x)))\geq ch^d>1.
$$
Taking $x= h{\bf n}$, where ${\bf n} \in {\Bbb Z}^d$, we see that all cubes of the form $h{\bf n} + [-h/2,h/2)^d$ contains at least one points of $\Lambda$, say $\lambda_{\bf n}$. Since $\Lambda$ is an F-spectrum,  $\{\lambda_{\bf n}\}_{{\bf n}\in{\Bbb Z}^d}$ is a Bessel set. By the stability under perturbation (see e.g., [DHSW, Proposition 2.3]) and
$$
|\lambda_{\bf n}-h{\bf n}|\ \leq \ {\rm diam}([-h/2,h/2)^d) \ = \ \sqrt{d}\ h,
$$
we conclude that  $h{\Bbb Z}^d$ is also a Bessel set of $L^2(\mu)$. As a Bessel set is invariant under translation, we see that the  finite union  ${\Bbb Z}^d = \bigcup_{{\bf k}\in \{0, \dots , h-1\}^d}(h{\Bbb Z}^d+{\bf k})$ is again a Bessel set of $L^2(\mu)$, which proves the claim.

 Now consider
$$
G(x): = \sum_{{\bf n}\in{\Bbb Z}^d} |\widehat{\mu}(x+{\bf n})|^2.
$$
$G$ is a periodic function (${\rm mod} \  {\Bbb Z}^d$). As ${\Bbb Z}^d$ is a Bessel set, applying the definition to  $e^{2\pi i \langle x, \cdot\rangle}$, we see that $G(x)\leq B<\infty$. Hence $G\in L^1([0,1)^d)$ and
$$
\int_{{\Bbb R}^d}|\widehat{\mu}(x)|^2dx = \sum_{{\bf n}\in{\Bbb Z}^d} \int_{[0,1)^d} |\widehat{\mu}(x+{\bf n})|^2dx =\int_{[0,1)^d}|G(x)|dx<\infty.
$$
This means that $\widehat{\mu}\in L^2({\Bbb R}^d)$, which implies $\mu$ must be absolutely continuous. This is a contradiction.

\bigskip

\noindent(iii) \  If $\mu$ is absolutely continuous, then the density function must be bounded above and below almost everywhere on the support of $\mu$ [L, Theorem 1.1]. Hence, $\Lambda$ is an F-spectrum of $L^2(\Omega)$, where $\Omega$ is the support of $\mu$. By Landau's density Theorem, $D^{-}\Lambda\geq{\mathcal L}(\Omega)>0$.
\eproof

\bigskip

Now it is easy to conclude that an F-spectral measure is of pure type.

\bigskip

\noindent{\bf Proof of Theorem \ref{th1.1}.} First let us assume that if $\mu$ is decomposed into non-trivial discrete and continuous parts, $\mu = \mu_d+\mu_c$. Let  $\Lambda$ be an F-spectrum of $\mu$. As $L^2(\mu_d)$ and $L^2(\mu_s)$ are non-trivial subspaces of $L^2(\mu)$, it is easy to see that $\Lambda$ is also an F-spectrum of both $L^2(\mu_d)$ and $L^2(\mu_c)$. Then $\#\Lambda<\infty$ by Proposition \ref{th2.1}(i); but $\#\Lambda=\infty$ since $L^2(\mu_c)$ is an infinite dimensional Hilbert space. This contradiction shows that $\mu$ is either  discrete or purely continuous.

\medskip

 Suppose $\mu$ is continuous and has non-trivial singular part $\mu_s$ and  absolutely continuous part $\mu_a$.  By applying the same argument as the above, $\Lambda$ is an F-spectrum of $L^2(\mu_s)$ and $L^2(\mu_a)$.  This is impossible in view of the Beurling density of $\Lambda$ in Proposition \ref{th2.1}(ii) and (iii).\eproof

\bigskip

The following corollary is immediate from Theorem \ref {th1.1}.

\medskip

\begin{coro}
A spectral measure or an R-spectral measure must be of pure type.
\end{coro}

\bigskip

\section {Discrete measures}
\setcounter{equation}{0}

In this section, we will show that all discrete measures on ${\Bbb R}^d$ are R-spectral measures. By Proposition \ref{th2.1}(i), we only need to consider measures with finite number of atoms. Let ${\mathcal C} = \{c_0,...,c_{n-1}\}\subset {\Bbb R}^d$ be a finite set and let
\begin{equation}\label{eq3.1}
\mu = \sum_{c\in{\mathcal C}}p_c\delta_c \ , \quad \hbox {with}\ \ p_i >0 , \  \sum_{c\in{\mathcal C}}p_c =1.
\end{equation}
For $\lambda\in{\Bbb R}^d$, we denote the vector
$
[e^{2\pi i \langle\lambda,\ c_0\rangle},...,e^{2\pi i \langle\lambda,\ c_{n-1}\rangle}]^t
$ by  ${\bf v}_\lambda$.

\bigskip

\begin{prop}\label{th3.1}
Let ${\mathcal C} =\{c_0,...,c_{n-1}\}\subset {\Bbb R}^d$ and let $\mu$ be as in (\ref {eq3.1}). Let $\Lambda = \{\lambda_0,...,\lambda_{m-1}\} \subset {\Bbb R}^d$ be another finite set. Then

\vspace{0.2cm}

\ {\rm (i)} \ $\Lambda$ is an F-spectrum of $\mu$ if and only if \ {\rm span}$\{{\bf v}_{\lambda_0},...,{\bf v}_{\lambda_{m-1}}\} = {\Bbb C}^n$.

\vspace {0.2cm}

{\rm (ii)} \ $\Lambda$ is an R-spectrum of $\mu$  if and only if $m=n$ in the above identity.
\end{prop}

\medskip

\noindent{\bf Proof.}
(i) Suppose first $\Lambda$ is an F-spectrum of $\mu$. Let ${\bf u}= [u_0, \cdots , u_{n-1}]^t$ be such that  $\langle{\bf u},{\bf v}_{\lambda_i}\rangle =0$ for all $i$.  Consider $f$ as a function defined on ${\mathcal C}$  with $f(c_i) =u_i/p_{c_i}$. By using the lower bound of the Fourier frame, we have
$$
A\sum_{c \in{\mathcal C}}|f(c)|^2p_c \ \leq \ \sum_{\lambda\in\Lambda}|\sum_{c\in{\mathcal C}}f(c)e^{2\pi i  \langle\lambda, c\rangle}p_c|^2
\ = \ \sum_{\lambda\in\Lambda} |\langle{\bf u}, {\bf v}_{\lambda}\rangle|^2 \ =\ 0.
$$
It follows that $f(c) =0$ for all $c\in{\mathcal C}$,   hence ${\bf u} = {\bf 0}$ and the necessity  in (i) follows

\medskip

Conversely, the assumption implies that  the vectors ${\bf v}_{\lambda_0},...,{\bf v}_{\lambda_{m-1}}$ form a frame on ${\Bbb C}^n$ (see [Chr, Corollary 1.1.3]),  i.e., there exist $A, B>0$ such that for all ${\bf u}=[u_0,...,u_{n-1}]^t\in{\Bbb C}^n$
$$
A\sum_{i=0}^{n-1}|u_i|^2\leq\sum_{\lambda\in\Lambda}|\langle{\bf u},{\bf v}_{\lambda}\rangle |^2\leq B\sum_{i=0}^{n-1}|u_i|^2.
$$
For any $f\in L^2(\mu)$, we take ${\bf u} = [f(c_0)p_{c_0},...,f(c_{n-1})p_{c_{n-1}}]^t$, we see that $\Lambda$ is a frame with lower bound \   $(\min_i {p_i}) A$ \  and upper bound  \ $(\max_i{p_i}) B$.

\medskip

(ii) is clear from (i).
\eproof

\bigskip

\noindent{\bf Proof of Theorem \ref{th1.2}.}  Let $\C =\{c_0, \cdots , c_{n-1}\}$. We first establish the theorem for  ${\mathcal C} \subset {\Bbb R}^1$. Let $W = {\rm span}\{{\bf v}_{\lambda}: \lambda\in{\Bbb R}^1\}$, it suffices to show that $W = {\Bbb C}^n$. Then we can select $\{\lambda_0,...,\lambda_{n-1}\}\subset {\Bbb R}^1$ so that $\{{\bf v}_{\lambda_0},...,{\bf v}_{\lambda_{n-1}}\}$ a basis of ${\Bbb C}^n$. The theorem for ${\Bbb R}^1$ will follow from Proposition \ref{th3.1}(ii).

\medskip

To see $W = {\Bbb C}^n$, it suffices to show that if $\langle{\bf u},{\bf v}_{\lambda}\rangle = 0$ for all $\lambda\in{\Bbb R}$, then ${\bf u} = {\bf 0}$. To this end, we write ${\bf u} = [u_0,...,u_{n-1}]^t$, and the given condition  is
\begin{equation*}\label{eq3.2}
\sum_{i=0}^{n-1}u_ie^{2\pi i\lambda c_i} = 0 .
\end{equation*}
We differentiate the expression with respect to $\lambda$ for $k$ times with $k = 1,...,n-1$, then
$$
\sum_{i=0}^{n-1} u_ic_i^ke^{2\pi i\lambda c_i} = 0.
$$
This means
$$
\left[
  \begin{array}{cccc}
    1 & 1 & \cdots & 1 \\
    c_0 & c_1 & \cdots & c_{n-1} \\
     \vdots& \vdots & \vdots & \vdots \\
    c_0^{n-1} & c_1^{n-1} & \cdots & c_{n-1}^{n-1} \\
  \end{array}
\right]\cdot\left[
              \begin{array}{c}
                u_0e^{2\pi i \lambda c_0} \\
                u_1e^{2\pi i \lambda c_1} \\
                \vdots \\
                u_{n-1}e^{2\pi i \lambda c_{n-1}} \\
              \end{array}
            \right] = {\bf 0}.
$$
As all $c_i$ are distinct, the Vandermonde matrix is invertible. Hence,
$$
[u_0e^{2\pi i \lambda c_0},u_1e^{2\pi i \lambda c_1},...,u_{n-1}e^{2\pi i \lambda c_{n-1}}]^t = {\bf 0},
$$
and thus ${\bf u} = {\bf 0}$. This completes the proof of the theorem for ${\Bbb R}^1$.

\vspace{0.2cm}

On ${\Bbb R}^d$, we note that by Proposition \ref{th3.1}, $\{c_0,...,c_{n-1}\}$ admits an R-spectrum $\{\lambda_0,...,\lambda_{n-1}\}$ if and only if $\{Qc_0,...,Qc_{n-1}\}$ admits an R-spectrum $\{Q\lambda_0,...,Q\lambda_{n-1}\}$, where $Q$ is any orthogonal transformation on ${\Bbb R}^d$. Now given any $\{c_0,...,c_{n-1}\}$ on ${\Bbb R}^d$, we let $\ell_{ij}$ be the line passes through two points $c_i$, $c_j$, and choose a line $\ell$ such that $\ell$ is {\it not}  perpendicular to any $\ell_{ij}$. Apply an orthogonal transformation $Q$ so that the first axis coincides with the direction of $\ell$. In this way the construction shows that the first coordinates of $Qc_0,...,Qc_{n-1}$ are all distinct.

We then apply the same argument as in ${\Bbb R}^1$ above, using partial differentiation with respect to the first coordinates which are all distinct, the Vandermonde matrix is invertible, hence the theorem follows.
\eproof

\bigskip

\noindent {\bf Remark 1}. \ If $c_0,...,c_{n-1}$ have rational coordinates,  then we can choose elements $\Lambda$ to have rational coordinates also. To see this, by multiplying an integer, we can assume that  $\{c_0,...,c_{n-1}\}$ are in ${\Bbb Z}^d$, we consider the determinant function
$$
\varphi(\lambda) = \varphi(\lambda_0,...,\lambda_{n-1})
= \det\left[
\begin{array}{ccc}
e^{2\pi i \langle \lambda_0, c_0\rangle} & \cdots & e^{2\pi i \langle \lambda_1, c_{n-1}\rangle} \\
& \ddots &  \\
 e^{2\pi i \langle \lambda_{n-1}, c_0\rangle}  & \cdots & e^{2\pi i \langle \lambda_{n-1}, c_{n-1}\rangle}  \\
 \end{array}
 \right]
$$
with $\lambda = (\lambda_0,...,\lambda_{n-1})$ on ${\Bbb R}^{dn}$.  Then $\varphi (\lambda)$  is a trigonometric polynomial on ${\Bbb R}^{dn}$, whose zero set is a closed set of Lebesgue measure zero. We can choose $\lambda$ so that $\varphi(\lambda)\neq 0 $ and $\lambda$ is rational, and Proposition \ref {th3.1}(ii) shows that $\Lambda = \{\lambda_0,...,\lambda_{n-1}\}$ is an $R$-spectrum will rational coordinates.

\bigskip

\noindent {\bf Remark 2}.  The R-spectrum shown in Theorem \ref{th1.2} is not explicit. It is also not easy to see whether a given set $\Lambda$ is an R-spectrum since the invertibility of the matrix is not easy to establish in general. A probabilistic approach of finding such $\Lambda$ in the case of trigonometric polynomials was given in [BG]. The work gave a theoretical background on the theory of reconstruction of multivariate trigonometric polynomials via random sampling sets.

\bigskip

To carry out Theorem 1.2 further, we consider the condition that a discrete measure to be an orthogonal spectral measure. We will restrict our consideration on the one dimensional case, and by translation, we can assume,  without loss of generality, that $\C \subset {\Bbb Z}^+$ and $0 \in \C$. The mask polynomial of $\mu = \sum_{c \in \C} p_c \delta_c$ is
$$
m_{\C,P} (x) \ = \ \widehat \mu (x) \  = \  {\sum}_{c \in \C} p_c e^{2 \pi i cx}.
$$
In case $P$ is a set of equal probability, then we just use the notation $m_\C (x)$. We call a set $\Lambda$  a \textit{bi-zero set} of $m_{\C, P}$  if  $0 \in \Lambda$ and  $m_{\C, P}( \lambda_i - \lambda_j)= 0$ for distinct $\lambda_i , \lambda_j \in \Lambda$. It is clear that such $E(\Lambda)$ is an orthogonal set in $L^2(\mu)$.

\bigskip

\begin{prop}\label{th3.4} Let ${\mathcal C} \subset {\Bbb Z}$ be a finite set,  and let \ $\mu  = \sum_{c \in {\mathcal C}} p_c \delta_c$. Then  $\mu$ is a spectral measure if and only if there is a bi-zero set $\Lambda$ of
 $m_{\C, P}$ and \ $\#\C = \#\Lambda$. In this case, all the $p_c$'s are equal.
\end{prop}

\medskip

\noindent{\bf Proof.} Note that $\mu$  is a spectral measure if and only if there exists a set $\Lambda  =\{\lambda_1, \ldots,
\lambda_n\}$ with $n=\#\C$ such that $\widehat{\mu}(\lambda_i-\lambda_j) = 0$ for all $i\neq j$. Since $\widehat{\mu}(x)=m_{\C, P}(x)$, this is equivalent to $\Lambda$ is a bi-zero set  of
 $m_{\C, P}$ and \ $\#\C = \#\Lambda$.

%Taking $x=-\lambda_j$ we obtain  $m_{\C, P}(\lambda_i-\lambda_j)=0$ for $i\ne j$. Hence, $\Lambda$ a bi-zero set of $m_{\C, P}$
%with $\#\Lambda=\#\C$. Conversely, if $\Lambda$ is a bi-zero set of $m_{\C, P}$
%with $\#\Lambda=\#\C$, then
% $\{e^{2\pi i \langle\lambda, x\rangle}\}_{\lambda \in \Lambda}$ is an orthogonal basis.

\medskip

To see that all the $p_c$ are equal, we put $f = \chi_{c}$ into the Parseval's identity.
$$
\sum_{\lambda\in\Lambda}|\langle f, e^{2\pi i \langle\lambda,x\rangle}\rangle|^2 = \|f\|^2.
$$
We obtain $\sum_{\lambda\in\Lambda}p_c^2 = p_c$. Hence, $p_c = 1/\#\Lambda = 1/\#\C$.
\eproof
%
% we note that $
%\sum_{j=1}^n|m_{\C, P}(x+\lambda_j)|^2\equiv 1 \ $ by Proposition \ref{th3.3}.
% We take $x =0$, then
%$$
%1\equiv\sum_{j=1}^n|\sum_{c\in\C}p_ce^{2\pi i \lambda_jc}|^2=
% \sum_{c,c'\in\C}p_cp_c'\sum_{j=1}^ne^{2\pi i \lambda_j(c-c')}.
%$$
%As $\Lambda$ is a bi-zero set of $m_{\A,P}$, the matrix $[\sqrt{p_c}e^{2\pi i \lambda c}]_{\lambda\in{ \Lambda},c\in\C}$ is unitary. Hence, the inner sum of the above expression is non-zero if and only
%if $c=c'$. Therefore
%$$
%\frac{1}{\#\C}=\sum_{c\in\C}p_c^2.
%$$
%By the Cauchy-Schwarz inequality, $1=(\sum_{c\in\C} p_c)^2\leq\sum_{c\in\C}p_c^2\cdot(\#\C)$, this means
% equality  holds, which implies $p_c = 1/\#\C$.

\bigskip

In the following we  use Proposition \ref{th3.4} to obtain explicit expressions of $\C$ with $\#\C \leq 4$ that are discrete spectral measures. It is trivial to check that when  $\#\C =1, 2$,  the associated $\mu$ is always a spectral measure.

\medskip

 \begin{example} \label{th3.5} {\it Let $\C=\{c_0=0, c_1, c_2\}\subset\Z^+$ with $\gcd(\C)=1$. Then
 $\mu = \sum_{c\in \C} \delta_c$ is a spectral measure if and only if $
 c_2\equiv2c_1 \ (\mbox{mod}\,\, 3)$ (i.e., $\C$ is a complete residue (mod 3)).}
 \end{example}

{\noindent \bf Proof.}  For the sufficiency, we write $c_1=3k+i,\ i = 1,2$ and $k \in \N$.  Then $c_2=3l+2i$ and
$$
\{0, c_1, c_2\}=\{0, 1, 2\} \ (\mbox{mod}\,\, 3).
$$
It is direct to check that  $\Lambda =\{0, \frac 13, \frac 23\}$  is a
bi-zero set of $m_{\C}(x)$ with $\#\Lambda=\#\C$and hence $\mu$ is a spectral measure.

\medskip
 For the necessity, we let  $\Lambda =\{0, \lambda_1, \lambda_2\}$ be such that
$m_\C(b_1)=m_\C(b_2)=m_\C(b_2-b_1)=0$. Note that
$m_\C(x)=1+e^{2\pi i c_1x}+e^{2 \pi i c_2x}$. Then $m_\C(x)$
has roots in $(0, 1)$ if and only if  $e^{2 \pi i c_1 x}=e^{{2\pi i}/3},
e^{2 \pi i c_2}=e^{{4\pi i}/3}$ (or the other way round). Hence there exists $k, l\in \Z^+$ such that
$$
 2\pi c_1x =2k\pi+\frac23\pi, \qquad 2\pi c_2x =2l\pi+\frac 43\pi.
$$
It follows that $x=\frac{3k+1}{3c_1}=\frac{3l+2}{3c_2}$.
Since $\gcd(c_1, c_2)=1$, we have
$3k+1 =c_1m$ and   $3l+2 =c_2m$.
Hence $3\nmid m$, and $3\mid (c_2-2c_1)$. This implies the sufficiency.
\eproof

\bigskip

\begin{example} \label{th3.6} {\it Let $\C=\{c_0=0, c_1, c_2, c_3\}\subset\Z^+$ with $\gcd(\C)=1$. Then
 $\mu$ is a spectral measure if and only if after rearrangement,
 $c_1$ is even, $c_2, c_3$ are odd, and $c_1 = 2^\alpha (2k+1)$, $c_2-c_3 =2^\alpha (2\ell+1)$ for some $\alpha > 0$.}
 \end{example}

\medskip

{\noindent \bf Proof.}  We first prove the necessity. The mask polynomial of $\mu$ is
$m_\C(x)=1+e^{2\pi i c_1 x}+e^{2 \pi i c_2x}+e^{2 \pi i c_3 x}$. That  $m_\C(x) =0$ implies
\begin {equation} \label {eq3.5}
|1+e^{2 \pi i c_1x}|=|1+e^{2 \pi i (c_3-c_2)x}|,
\end{equation}
which yields  (i) $e^{2 \pi i c_1x}=e^{2 \pi i (c_3-c_2)x}$ or (ii)
$e^{2 \pi i c_1x } = e^{- 2\pi i (c_3-c_2)x}$. Putting (i) into $m_\C(x) = 0$, we have $(1+ e^{2\pi i c_1x})(1 +e^{2\pi i c_2x}) =0$. Hence we have two sets of equations:
\begin{equation}\label{eq3.6}
2c_1x = 2k+1; \quad  2(c_3-c_2)x = 2l+1.
\end{equation}
or
\begin{equation}\label {eq3.7}
2c_2x = 2k+1; \quad  2(c_3-c_1)x = 2l+1.
\end{equation}
From (\ref{eq3.6}), we have
$
x=\frac {2k+1}{2c_1}= \frac {2l+1}{2(c_3-c_2)}.
$
Let $a = \gcd(c_1, c_2-c_3)$. It is easy to show that there exists $m$ such that
\begin {equation} \label {eq3.8}
2k+1 = \frac {c_1m}a \ ,  \qquad  2l+1 =\frac {(c_3-c_2)m} a \ .
\end{equation}
Hence $m, \ {c_1}/a,\ {(c_3-c_2)}/a$ \ must be odd. Also note that $\gcd ( c_1, c_2, c_3) =1$,  it follows from a direct check of the above that two of the $c_1, c_2, c_3$ must be odd, and one must be even (all three cases can happen).

The same argument applies to (\ref {eq3.7}) and to (ii). The last statement also follows in the proof.

\medskip

To prove the sufficiency, we first observe from the above that for $c_1$ even, $c_2, c_3$ odd, there are solutions  $x_1, x_2 \in (0,1)$ from (i) (see (\ref {eq3.6})-(\ref {eq3.8})):
$$
x_1 = \frac{2i+1}{2a}, \ \ 0\leq i < a \ ;  \qquad
x_2 = \frac{2j+1}{2b}, \ \ 0 \leq j < b \ ,
$$
where $\gcd(c_1, c_2-c_3) $ as above,  and $ b= \gcd(c_2, c_3-c_1)$.  Since $b$ is odd, we can take $2j+1 =b$, so that $x_2=\frac 12$ is a solution of $m_\C(x)=0$. Let
$$
 \lambda_1 = \frac12, \quad  \lambda_2=\frac 1{2a}, \quad  \lambda_3=\frac {2^\alpha\gcd(r,s)+1}{2a} \ ,
$$
where $c_1 = 2^\alpha r, c_2 = 2^\alpha s$ and $r, s$ are odd integers as in the assumption. We claim that $\Lambda = \{0, \lambda_1, \lambda_2, \lambda_3\}$ is a bi-zero set of
$m_\C(x)$. Indeed, since $a=2^\alpha\gcd(r, s)$,
$
\lambda_1-\lambda_2=(a-1)/{2a}
$
is of the form $x_1$ for  $i=2^{\alpha-1}\gcd(r, s)-1$, $\lambda_3-\lambda_1=\lambda_2$ and
$\lambda_3-\lambda_2=\lambda_1$, the claim follows.\eproof

\bigskip

For $\#\C$ large, it is difficult to evaluate the zero set of the mask polynomial. However there is a number-theoretical approach to study such zeros related to the spectrum and integer tiling,  this is to be discussed in the next section and in the Appendix.

\bigskip
\bigskip

\section { A connection with integer tiles}
\setcounter{equation}{0}

\bigskip

In this section we will give a brief discussion of the relationship between discrete spectral measures and integer tiles, and provide the tools we need in the next section.  Let $\A \subset {\Bbb Z}^+$ and assume that $0 \in \A$, we say that $\A$ is an {\it integer tile} if there exists $\mathcal T$ such that  $\A \oplus {\mathcal T} = {\Bbb Z}$, i.e.,  $\A + [0,1]$ tiles ${\Bbb R}$.  Equivalently, $\A$ is a tile if there exists ${\B}$ and $n$ such that
\begin {equation} \label {eq4.1}
\A \oplus \B \equiv {\Bbb Z}_n \  ({\rm mod } \ n).
\end{equation}

Recall that the Fuglede conjecture asserts that for $\Omega \subset {\Bbb R}^d$ with positive measure, $\Omega$ is a translational tile if and only if the restriction of the Lebesgue measure $\mathcal L |_\Omega$ is a spectral measure. Although the conjecture is proved to be false in either direction ([T], [KM]), it remains unanswered for dimension $1$ and $2$, and for some special classes of tiles in any dimension.

\bigskip

Let $\A$ be a finite subset in ${\Bbb Z}$,  then the Fuglede conjecture reduces to {\it $\A$ is an integer tile if and only if ${\mathcal A}+[0,1]$ is a spectral set, i.e., $\mathcal L|_{\A + [0,1]}$ is a spectral measure}. It is also known that the latter part is also equivalent $\delta_\A = \sum_{a\in \A} \delta_a$ is a discrete spectral measure  [LW]. This also follows from Theorem \ref{th5.3} in Section 5. In Example \ref {th3.5}, the spectral condition for $\# \C =3$ is equivalent to ${\mathcal C}$ is a complete residue (mod $3$), which trivially satisfies (\ref {eq4.1}). Hence the conjecture is true for $\# \C =3$.  In Example \ref{th3.6}, the spectral condition is equivalent to
$$
\C = \{0,\ 2^{\alpha}(2k+1),\ 2r+1,\ 2r+1+2^{\alpha}(2\ell+1) \}
$$
for some non-negative integers $k,r,\ell$.
If we let
$
{\mathcal B} = \{0,2\}\oplus...\oplus \{0,2^{\alpha-1}\},
$
then it is direct to check that ${\mathcal C}\oplus{\mathcal B} \equiv{\Bbb Z}_{2^{\alpha+1}}$  ({\rm mod } \ $2^{\alpha+1}$). Hence by (\ref {eq4.1}), the conjecture is true for $\#\C =4$. Actually, by using some deeper number-theoretic argument (see Appendix), it can be shown that if $\# \C = p^\alpha q^\beta$ where $p,q$ are distinct primes, then ${\mathcal C}$ is an integer tile implies it is a spectral set ([CM], [{\L}a]).

\bigskip

 The following is a useful sufficient condition of a discrete spectral measure. The condition trivially imply the underlying set is an integer tile.

\bigskip

\begin{theo} \label{th4.1} \ Let $\A \subset {\Bbb Z}^+$ be a finite set with $0\in \A$. Suppose there exists $\B \subset {\Bbb Z}^+$ such that
$$
\A \oplus \B = {\Bbb N}_n
 $$
 where ${\Bbb N}_n = \{0, \cdots , n-1\}$. Then the discrete measure $\delta_{{\mathcal A}} = \sum_{a\in \A} \delta_a$ (with equal weight) is a spectral measure with a spectrum contained in $\frac{1}{n}{\Bbb Z}$.
\end{theo}
\medskip

The theorem was due to [PW] (and also in [DJ]), and the proof  involves an inductive construction of the spectrum. The spectrum is implicit and the proof is long. We will provide an alternative proof using the properties of the root of unity as the zeros of the mask polynomial. The framework is from [CM]  and the spectrum is explicitly given in [{\L}a].  Because of the number-theoretical notations and techniques,  we will leave the details in the Appendix.

\bigskip

Finally, we state a related theorem of the self-similar measures which follows from the known results, and will be needed in the next section.

\bigskip

\begin{theo} \label{th4.2} \ Let $\A \subset {\Bbb Z}^+$ be a finite set with $0\in \A$. Suppose there exists $\B \subset {\Bbb Z}^+$ such that $\A \oplus \B = {\Bbb N}_n$.   Let $\mu$ be the the self-similar measure satisfying
$$
\mu( \cdot ) = \frac 1{\# A} \sum_{a \in \A} \mu (n\cdot -a).
$$
Then $\mu$ is a spectral measure. Moreover, if  gcd $(\A)=1$, then the spectrum $\Lambda$ of $\mu$  can be chosen to be in ${\Bbb Z}$.
\end{theo}

\medskip

\noindent{\bf Proof.}  Denote the spectrum of $\delta_{\A}$ in Theorem \ref{th4.1} by ${\mathcal S}$ with ${\mathcal S}\subset \frac 1n{\Bbb Z}$, then  $(\A,{\mathcal S})$ form a \textit{compatible pair} as in [{\L}aW1], i.e.,
$$
\frac{1}{\sqrt{\#\A}}[e^{2\pi i as}]_{a\in{\mathcal A}, s\in{\mathcal S}}
$$
is a unitary matrix.  The theorem  follows from [{\L}aW1, Theorem 1.2]. For the last part, since we can change the residue representatives of $\Gamma$ with ${\mathcal S}  = \frac{1}{n}\Gamma$ and $\Gamma \subset\{-(n-2),...,n-2\}$. With gcd $(\A)=1$,  Theorem 1.2 in [{\L}aW1] states that
$$
\Lambda = \Gamma\oplus n\Gamma\oplus n^2\Gamma\oplus....
$$
  is a spectrum. This spectrum clearly lies in ${\Bbb Z}$. \eproof

\bigskip

\noindent {\bf Remark}. For the $1/4$-Cantor measure, $\mu(\cdot) = \frac 12 \mu( 4\cdot ) + \frac 12 \mu (4\cdot - 2)$.  It is easy to compute the Fourier transform is
$
\widehat{\mu}(\xi) = e^{2\pi i \frac{1}{3}\xi}\prod_{j=1}^{\infty}\cos({2\pi \xi}/{4^j})$ and the zero set of $\widehat \mu$ is ${\mathcal Z}_{\mu} = \{4^{j}a: a \ {\mbox{is odd and}} \  j\geq0\}$.
Note that $\A = \{0,2\}$ and the condition of the theorem is satisfied,  the spectrum of $\mu$ can be taken as [JP]
$$
\Lambda = \{0,1\}\oplus 4\{0,1\}\oplus4^2\{0,1\}\oplus \cdots \subset {\mathcal Z}_\mu.
$$

\medskip

However the condition ${\mathcal A}\oplus{\mathcal B} ={\Bbb N}_n$ in Theorem \ref{th4.2} is quite restrictive. For the $1/6$-Cantor measure, $\mu(\cdot) = \frac 12 \mu( 6\cdot ) + \frac 12 \mu (6\cdot - 2)$, according to [JP],  it is again a spectral measure and the spectrum is
 $$
\Lambda = \frac 32 \big (\{0,1\}\oplus 6\{0,1\}\oplus6^2\{0,1\}\oplus \cdots \big ).
$$
But for $\mathcal A = \{0, 2\}$ in this case, we cannot find $\B$ so that $\mathcal A \oplus \mathcal B = {\Bbb N}_6$. Also, $\mu$ does not admit spectrum $\Lambda ' \subset {\Bbb Z}$.  The proof is as follows: If so, observe that
$$
\Lambda ' \ \subset \  {\mathcal Z}_{\mu} \ = \ \{6^{j}a/4: a \ {\mbox{is odd and}} \  j\geq1\}.
$$
 As $\lambda$ is an integer, we see that for  $\lambda\in\Lambda'$, $\lambda = {6^na}/{4}$, and $n \geq 2$ necessarily.  Let $x = {3}/{2}$, then \ $x\in{\mathcal Z}_{\mu} \setminus \Lambda'$ \ and \ $x-\lambda = {6(1-6^{n-1}a)}/{4}\in {\mathcal Z}_{\mu}$. This means $\sum_{\lambda\in\Lambda}|\widehat{\mu}(x-\lambda)|^2 =0$, which shows that $\Lambda'$ cannot be a spectrum by Proposition \ref{th3.3}.

\bigskip
\bigskip

\section {Convolutions}
\setcounter{equation}{0}

\medskip

Let $\nu$ be a probability measure with compact support $\Omega\subset[0,1]$ and let $\eta$ be a discrete probability measure with support on $\C \subset {\Bbb Z}^+$ and probability weight $P$, i.e. $\eta = \delta_{\C, P} = \sum_{c \in \C} p_c\delta_c$.  Then $\mu=\eta\ast \nu$ has support on $\C+ \Omega$.  Given a non-negative integer $q$, we let $\eta_ q= \delta_{q\C,P}$.
%and
%satisfies
%\begin{equation}\label{eq5.1}
%\mu(E)=\eta\ast \nu(E)= \sum_{c\in \C}p_c \nu(E-c)
%\end{equation}
%for any Borel set $E$.

%\bigskip
%
%\begin{lem}\label{th5.1}
%Let ${\mathcal A}\subset {\Bbb N}$ be a finite set. Then there
%exists ${\mathcal S}$ such that $\#{\mathcal S}=\#{\mathcal A}$ and
%$$
%\det[e^{2\pi i ab}]_{a\in{\mathcal A},s\in{\mathcal S}}\ \neq \ 0.
%$$
%Moreover, we can choose ${\mathcal S}$ so that every element is
%rational.
%\end{lem}

\bigskip

\begin{theo}\label{th5.2} Let  $\nu$ be an R-spectral measure with a spectrum $\Gamma$ and assume that there exists an integer $q\geq 1$ such that $q\Gamma \subseteq \Z$.  Then  $\mu := \eta_q \ast \nu$ is an R-spectral measure.
\end{theo}

\bigskip

 \noindent{\bf Proof.}  We write  $ \A = q\C= \{0 = a_0, a_1, \ldots, a_{k-1}\}$. By Theorem \ref{th1.2}, there exists an  R-spectrum of $\A$ which we  denote it as ${\mathcal S} =\{ 0=s_0,s_1,\ldots, s_{k-1}\}$. By Proposition \ref{th3.1}(ii), we see that
$$
\det[e^{2\pi i as}]_{a\in{\mathcal A},s\in{\mathcal S}}\ \neq \ 0.
$$
We will show that
${\mathcal S} \oplus \Gamma$ is an R-spectrum of $\mu$.

\vspace {0.15cm}

Since $\Omega = {\rm supp} (\nu)\subset [0,1]$, for any $f\in L^2(\mu)$, $f$ is uniquely determined by the vector-valued
function $[f(x+a_0), \cdots, f(x+a_{k-1})]^t$ on  $\Omega$ . Let
$M=\big[e^{2\pi i a s}\big]_{a\in\A, s\in{\mathcal S}}$, it is invertible.  We define
$$
[g_0(x), \cdots, g_{k-1}(x)]^t = M^{-1}[f(x+a_0), \cdots,
f(x+a_{k-1})]^t\ , \quad x \in \Omega\ .
$$
Clearly $g_j\in L^2(\nu)$ for $0\le j\le k-1$.  It is easy to see $s_j+\Gamma$ is also an R-spectrum of $\nu$, so that $g_j$ can be uniquely expressed as
$$
g_j(x)=\sum_{\gamma\in\Gamma}c_{s_j+\gamma}e^{2\pi i
(s_j+\gamma)x}.
$$
Hence, we have
$$
M[g_0(x), \cdots, g_{k-1}(x)]^t=\left[\sum_{j=0}^{k-1}e^{2\pi i
a_0s_j}g_j(x), \cdots,\sum_{j=0}^{k-1}e^{2\pi i
a_{k-1}s_j}g_j(x)\right]^t,
$$
and therefore
\begin{equation} \label{eq5.2}
f(x+a_i)=\sum_{j=0}^{k-1}e^{2\pi i a_is_j}\sum_{\gamma\in\Gamma}c_{s_j+\gamma}e^{2\pi i
(s_j+\gamma)x}=\sum_{j=0}^{k-1}\sum_{\gamma\in\Gamma}c_{s_j+\gamma}e^{2\pi
i (s_j+\gamma)(x+a_i)}.
\end{equation}
 Note that the last equality follows from $\gamma a_i = (q\gamma) c_j$ is an integer by the assumption $q\Gamma\subset{\Bbb Z}$.  By a change of variable with $y = x+a_i$ for each $i$, we have
$$
f(y)=\sum_{j=0}^{k-1}\sum_{\gamma\in\Gamma}c_{s_j+\gamma}e^{2\pi
i (s_j+\gamma)y}, \qquad y\in  \C + \Omega = {\rm supp} (\mu).
$$
It is easy to see that the above representation  is unique, this
 means $E({\mathcal S}+\Gamma)$ is both a basis and a frame of $L^2(\mu)$. Hence $E({\mathcal S}+\Gamma)$ is a Riesz basis.
 \eproof

 \bigskip

 We now recall a general criterion of spectral measures due to Jorgensen and Pedersen [JP].

\medskip

\begin{prop}\label{th3.3} Let $\mu$ be a probability measure on ${\Bbb R}^d$ with compact support. Then $\Lambda$ is an orthogonal spectrum of $\mu$ if and only if
$$
Q(x)=\sum_{\lambda\in\Lambda}|\widehat{\mu}(x+\lambda)|^2\equiv 1 \ ,  \quad
 x\in \R.
$$
\end{prop}

\medskip

In particular, if $\mu = \sum_{c\in\C}p_c\delta_c$ is a discrete spectral measure with spectrum $\Lambda$, then $p_c = {1}/{\#\C}$ by Proposition \ref{th3.4} and
 $$
\sum_{\lambda\in\Lambda}|m_{\C,P}(x+\lambda)|^2\equiv 1. 
$$
To determine whether  $\mu$ in Theorem \ref{th5.2} is a spectral measure, we have the following simple characterization.

 \medskip

 \begin{theo}\label{th5.3} Let  $\nu$ be an R-spectral measure and suppose  $q {\mathcal Z}_\nu \subset {\Bbb Z}$. Then $\mu = \eta_q \ast \nu$ is a spectral measure if and only if both $\eta$  and $\nu$ are  spectral measures.
 \end{theo}

 \bigskip

\noindent {\bf Proof}. It is clear that $\eta$ is a spectral measure if and only if $\eta_q$ is also a spectral measure. We first prove the sufficiency. Let $\A=q\C$, and let
${\mathcal S} =\{0, s_1, \ldots, s_{k-1}\} $ be a bi-zero set of $m_{\A,
P}$ (Note that $P$ is a set of equal weights by Proposition \ref{th3.4}).  Let $\Gamma$ be a  spectrum of $\nu$, then $q\Gamma\subseteq
\Z$ by the hypothesis that $q {\mathcal Z}_\nu \subset {\Bbb Z}$. The  Fourier transform  of $\mu$ satisfies
$$
\widehat{\mu}(\xi)=m_{\A,\,P }(\xi)\widehat{\nu}(\xi),
$$
 By the spectral property of  ${\mathcal S}$  and $\Gamma$,
\begin{eqnarray*}
\sum_{0\le j\le
k-1}\sum_{\gamma\in\Gamma}|\widehat{\mu}(x+s_j+\gamma)|^2
&=&\sum_{0\le j\le k-1}\sum_{\gamma\in\Gamma}|m_{\A,\, P}(x+s_j+\gamma)|^2|\widehat{\nu}(x+s_j+\gamma)|^2\\
&=&\sum_{0\leq j\le k-1}\sum_{\gamma\in\Gamma}|m_{\A,\, P}(x+s_j)|^2|\widehat{\nu}(x+s_j+\gamma)|^2\\
&=&\sum_{0\le j\le k-1}|m_{\A,\, P}(x + s_j)|^2=1.
\end{eqnarray*}
Hence ${\mathcal S} \oplus \Gamma$ is an orthogonal spectrum of $\mu$ by Proposition \ref{th3.3}.

\vspace {0.2cm}

Conversely, suppose that $\Lambda$ is a spectrum of $\mu$ and without loss of generality assume $0\in \Lambda$.   Denote $x=\{x\}+ [x]$ where $[x]$ is the maximum integer which is less than or
equal to $x$. We claim that ${\mathcal S} =\{q^{-1}\{q\lambda\}:
\lambda\in \Lambda\}$ is a bi-zero set of $m_{\A, \, P}$. Indeed,  by writing \ $\lambda = q^{-1} \{q \lambda\}+q^{-1} [q\lambda]$, we have
\begin{equation} \label{eq5.3}
0 = \widehat{\mu}(\lambda) \ = \ m_{\A, \, P}\big (q^{-1}\{q\lambda\} + q^{-1} [q \lambda] \big )\widehat{\nu}(\lambda)
\ =\  m_{\A, \, P}(q^{-1}\{q\lambda\})\widehat{\nu}(\lambda)
\end{equation}
for each $\lambda \in\Lambda$.  Note that $\widehat{\nu}(\lambda)=0$ implies
$q \lambda\in \Z$ (by the assumption $q{\mathcal Z}_\nu \subset {\Bbb Z}$), so that $\{q\lambda\}=0$,
 (\ref {eq5.3}) implies that either $q^{-1} \{q\lambda\} =0$ or it is a root of $m_{\A, \, P}$. For any given distinct $q^{-1}\{q\lambda_1\}$, $q^{-1}\{q\lambda_2\} \in\mathcal S$ and $\{q\lambda_1\}> \{q\lambda_2\}$, we have $q^{-1}\{q(\lambda_1-\lambda_2)\}=q^{-1}(\{q\lambda_1\}-\{q\lambda_2\})$ is a root of $m_{\A, \, P}$. This proves the claim.

\vspace {0.2cm}

Let us write
$\Lambda=\bigcup_{j=0}^{k-1}(s_j+\Lambda_j), s_j \in {\mathcal S}$, where $\Lambda_j = \{q^{-1}[q\lambda]:q^{-1}\{q\lambda\} =s_j\}$. Since $\Lambda$ is a spectrum of $\mu$, we must have for all $\lambda_1,\lambda_2\in\Lambda_j$,
$$
0=\widehat{\mu}(\lambda_1-\lambda_2) = m_{\A, \, P}(\lambda_1-\lambda_2)\nu(\lambda_1-\lambda_2).
$$
But $a(q^{-1}[q\lambda]) \in{\Bbb Z} $ for all $a\in{\mathcal A} = q{\mathcal C}$, this shows $m_{\A, \, P}(\lambda_1-\lambda_2)\neq 0$ and hence $\nu(\lambda_1-\lambda_2) =0$. Therefore $E(\Lambda_j)$, $0\le j\le k-1$, are the orthogonal set of $\nu$. By the Bessel inequality, $\sum_{\lambda\in\Lambda_j}|\widehat{\nu}(x+s_j+\lambda)|^2\leq1$. Note further that ${\mathcal S} $ is a bi-zero set of $m_{\A,P}$. By Proposition \ref{th3.3}, we have
 \begin{eqnarray*}
 1&\equiv&
 \sum_{\lambda\in\Lambda}|\widehat{\mu}(x+\lambda)|^2=\sum_{j=0}^{k-1}
 \sum_{\lambda\in\Lambda_i}|\widehat{\mu}
 (x+s_j+\lambda)|^2\\
 &=&\sum_{j=0}^{k-1}\sum_{\lambda\in\Lambda_j}|m_{\A,\,
 P}(x+s_j+\lambda)\widehat{\nu}(x+s_j+\lambda)|^2\\
 &=&\sum_{j=0}^{k-1}\sum_{\lambda\in\Lambda_j}|m_{\A,\,
 P}(x+s_j)\widehat{\nu}(x+s_j+\lambda)|^2  \qquad   {( \mbox{since} \ a\lambda\in{\Bbb Z})}\\
 &\le& \sum_{j=0}^{k-1}|m_{\A,\,
 P}(x+s_j)|^2\leq 1.
 \end{eqnarray*}
 Hence ${\mathcal S}$ is the orthogonal spectrum of $\eta_q$ by Proposition \ref{th3.3} again, so that $\eta$ is a spectral measure. From the third line of the above, we also have
 $$
 1\equiv\sum_{j=0}^{k-1}|m_{\A,\,
 P}(x+s_j)|^2\sum_{\lambda\in\Lambda_j}|\widehat{\nu}(x+s_j+\lambda)|^2.
 $$
 With $\sum_{j=0}^{k-1}|m_{\A,\,
 P}(x+s_j)|^2 \equiv1$, we must have $\sum_{\lambda\in\Lambda_j}|\widehat{\nu}(x+s_j+\lambda)|^2\equiv1$. Hence, $\nu$ is a spectral measure and any one of the $\Lambda_i$ is a spectrum of $\nu$.
 \eproof

\bigskip

It has been an open question whether the $1/3$-Cantor measure has an F-spectrum (or even an R-spectrum). To a less extend, we do not know a non-trivial  singularly continuous R-spectral measure. In the following, we can make use of Theorems \ref{th5.2} and \ref{th5.3} to construct such measures.

\bigskip

\begin {example} \label{eq5.5}\ {\it  There exists a singularly continuous R-spectral measure which is not a spectral measure.}
\end{example}

\bigskip

\noindent {\bf Proof.} \  Consider the self-similar measure $\nu_{{\A},n}$ in Theorem \ref{th4.2} with ${\mathcal A}$ satisfying ${\mathcal A}\oplus{\mathcal B} = {\Bbb N}_n$ and gcd$({\mathcal A}) =1$.  It is a spectral measure and has a spectrum $\Gamma \subset {\Bbb Z}$. Moreover we claim that ${\mathcal Z}_\nu \subset {\Bbb Z}$. Indeed, observe that
$$
\widehat{\nu}(\xi) = \prod_{j=1}^{\infty}m_{\A}(\frac{\xi}{n^j}).
$$
where  $m_{\A}$ stands for the mask polynomial of $\A$ under equal weight. As  ${\mathcal A}\oplus{\mathcal B} = {\Bbb N}_n$, we have
$$
m_{\A}({\xi})m_{\B}(\xi) = 1+e^{2\pi i \xi}+ ... +e^{2\pi i (n-1)\xi}.
$$
The zero set of $m_{\A}$ on $[0,1)$ is a finite subset $Z\subset\{1/n, ..., n-1/n\}$. Let $Z' = Z+{\Bbb Z}$. This shows that ${\mathcal Z}_{\nu} = \bigcup_{j=1}^{\infty}n^jZ'$. This proves the claim and  the condition in Theorem  \ref{th5.3} holds (taking $q =1$).

Now we let $\eta = \delta_{\C, P}$ be a discrete measure with any finite set $\C$  of non-negative integers  and {\it non-uniform} weight $P$.  $ \mu = \eta \ast \nu_{{\A},n}$ is an R-spectral measure but not a spectral measure by Theorem \ref{th5.2} and Theorem \ref{th5.3}. These measures is clearly singular if $\#\A<n$. \eproof

\bigskip

Finally, if $E$ is a Borel set with positive Lebesgue measure, we use $L^2(E)$ to denote the square integrable functions on $E$.  We remark that $L^2(E)$ always have an F-spectrum, and the existence of orthogonal spectrum  is related to the  translational tile as in Fuglede's conjecture.
For R-spectrum, it is not known whether \textit{every Borel set $E$ with positive Lebesgue measure,  $L^2(E)$ has an  R-spectrum.} In regard to this we have the following simple result.

\medskip

\begin{coro}  If $E$ be an finite union of closed intervals with rational endpoints. Then $L^2(E)$ admits an R-spectrum.
\end{coro}

\medskip

\proof By the hypothesis and by suitably rescaling and translation,  there exist two integers $r$ and $s$ such that
$$
rE+s \ = \ [0, 1]+\A : = F,
$$
where  $0\in\A$ and $\A\in \Z^+$ is a finite set. By Theorem \ref{th1.3} with $\nu$ being the  Lebesgue measure on [0,1],  we see that $F$ has an R-spectrum, which implies  $L^2(E)$ also has an R-spectrum.  \eproof

\bigskip

We remark that similar results were obtained in [LS] who considered the problem from the sampling point of view and used techniques in complex analysis. We do not know whether the condition of rational endpoints can be removed. In [Le], the case when the end-points lying in certain groups was considered, and the above also follows as a corollary.

\bigskip
\bigskip

\section{Appendix: Proof of Theorem \ref{th4.1}}
\setcounter{equation}{0}

In this section, we will  prove Theorem \ref{th4.1} using a number-theoretic method. The setup is in [CM]. For the trigonometric polynomial  $m_{\mathcal A}(\xi) =\sum_{a\in{\mathcal A}}e^{2\pi i a\xi}$, it is convenient to replace by the polynomial $P_{\mathcal A}(x) =\sum_{a\in{\mathcal A}}x^a$. Recall that a cyclotomic polynomial $\Phi_n(x)$ is the minimal polynomial of the $n$-th root of unity. It follows that  $\Phi_s(x) | P_{\mathcal{A}}(x)$ is equivalent to $m_\A (s^{-1}) = 0$.

\medskip

For a finite set ${\mathcal A} \subset {\Bbb Z}^+$, we write $p$ as primes and define
\begin{equation}\label{eq6.1}
S_{\mathcal{A}} = \{ p^\alpha > 1 : \  \Phi_{p^\alpha}(x) | P_{\mathcal{A}}(x) \}
\quad \hbox {and} \quad
\widetilde S_{\mathcal A} = \{s>1: \ \Phi_s(x) | P_{\mathcal{A}}(x) \}
\end{equation}
and the following two conditions

\vspace {0.2cm}

 { (T1)} \  { \it
$\#{\mathcal{A}}=\prod_{s\in S_{\mathcal{A}}} \Phi_s(1)$,}

\vspace {0.1cm}

 { (T2)} \ {\it For any distinct prime powers $s_1, \ldots,
s_n \in S_{\mathcal A} $ , then $s_1 \cdots s_n \in {\widetilde
S}_{\mathcal{A}}$. }

\medskip

  Intuitively, $(T1)$ means if $\#{\mathcal A} = \prod_{j=1}^{k}p_j^{\alpha_j}$, then there is exactly $\alpha_j$ prime powers of $p_j$ in  $S_{\mathcal{A}}$ (since $\Phi_{p^r}(1) = p$). Also, $(T2)$ is a generalization from the basic identity $1+x+...+x^{n-1} = \prod_{d|n,d>1}\Phi_{d}(x)$ to $P_{{\mathcal A}}(x)$. It is known that if a set ${\mathcal A}$ satisfies these conditions, then it it is an integer tile, i.e., it  tiles ${\Bbb Z}$. Conversely, if a set ${\mathcal A}$ tiles ${\Bbb Z}$, then it satisfies $(T1)$;  for $(T2)$ it holds when $\#\A = p^\alpha q^\beta$ for $p, q$ distinct prime numbers, and is still an open question without the additional condition. Theorem \ref {th6.3} we are going to prove is a special case of this. In another direction {\L}aba {\rm [{\L}a]} showed that $(T1)$ and $(T2)$ imply the spectral property:

  \bigskip

\begin{prop}\label{th6.1}
 Suppose that a finite set of non-negative integers ${\mathcal
A}$ satisfies $(T1)$ and $(T2)$. Then $\delta_{\A}$ is a spectral measure and has a spectrum
$$
 {\mathcal S} = \{\sum_{s\in S_{{\mathcal A}}}{k_s}{s^{-1}}:\
k_s\in\{0,1,...,p-1\} \ if  \ s = p^{\alpha} \}.
$$
\end{prop}

\bigskip

The following are basic manipulation rules for the cyclotomic polynomials.

\medskip

\begin{lem}\label{th6.2} Let $p$ be a prime, then

\vspace {0.15cm}

(i) \ $\Phi_{s}(x^p) = \Phi_{sp}(x)$ if  $p|s$, and

\vspace {0.15cm}
(ii) $\Phi_{s}(x^p) = \Phi_{s}(x)\Phi_{sp}(x)$ if $p$ is not a factor of $s$;

%\vspace {0.15cm}
%\noindent  Moreover $\Phi_{p^k}(1) = p$,  and $\Phi_n(1) = 1$ if $n>1$ is not a prime power.
\end{lem}

\bigskip

It follows that if $p$ is a prime  and $\Phi_{p^{\alpha}}(x)$ divides  $P_{\A}(x^m)$ for some $m>0$, then  by Lemma \ref{th6.2}, we have
$\Phi_{p^{\alpha+1}}(x)$ divides $P_{\A}(x^{m/p})$ if $p|m$ (by (i)),  and
$\Phi_{p^{\alpha}}(x)$ divides $P_{{\mathcal A}}(x^m)$ if $p\nmid m$ (by (ii)).
 If we let  $\gamma = \max\{j: p^{j}|m\}$ (by convention $\gamma=0$ if $p\nmid m$) and repeat the above argument, we have
$
\Phi_{p^{\alpha+\gamma}}(x)| P_{{\mathcal A}}(x).
$
That is
\begin{equation} \label {eq6.5}
S_{m\A} = \{ p^{\alpha+\gamma}: \ p^\alpha \in S_\A ,  \ \gamma = \max \{j: p^j|m\}\}.
\end{equation}

\bigskip

 We reformulate Theorem \ref{th4.1} as the following

\medskip

\begin{theo} \label{th6.3}  Let ${\mathcal A}\oplus {\mathcal B}={\mathbb N}_n$, then ${\mathcal A}$ and
${\mathcal B}$ satisfies $(T1)$ and $(T2)$. Hence $\delta_{\mathcal A}$ is a spectral measure with a spectrum in $\frac{1}{n}{\mathbb Z}$.
\end{theo}

\medskip
\noindent {\bf Proof.} It is easy to check that both ${\mathcal A}$ and ${\mathcal B}$ satisfy condition $(T1)$ [CM]. It
remains to prove condition $(T2)$ holds. Without loss of generality we assume that
$1\in \A$. We use induction on the number of primes of $n$.
If $n$ consists only of one prime. Then ${\mathcal A}$  equals
$\{0,1,...,n-1\}$ and $(T2)$ holds trivially.

\medskip

  Let ${\mathcal A} \oplus {\mathcal B} = {\Bbb N}_n$, then by [Lo, Lemma 1],
there exists  $m\geq2$ with $m|n$ and ${\mathcal A}'\oplus{\mathcal
B}' = {\mathbb N}_{n/m}$ such that
  \begin{equation}\label{eq6.2}
  {\mathcal A} = m{\mathcal A}'+{\mathbb N}_m \quad \mbox{and} \quad {\mathcal B} = m{\mathcal B}'.
  \end{equation}
Induction hypothesis implies that $(T2)$ holds for ${\mathcal A}'$ and ${\mathcal B}'$.

\medskip

To show that ${\mathcal B}$ has $(T2)$, we observe that  $P_{{\mathcal B}}(x)
=P_{{\mathcal B}'}(x^m)$. By (\ref {eq6.5})
$$
 S_{{\mathcal B}} = \{q^{\beta+\gamma}:q^{\beta}\in S_{{\mathcal B}'},  \ \gamma =\max\{j:q^{j}|m\} \}.
$$
 Taking distinct prime powers in $S_{{\mathcal B}}$, say $q_1^{\beta_1+\gamma_1},..., q_{\ell}^{\beta_{\ell}+\gamma_{\ell}}$. By $(T2)$
 of ${\mathcal B}'$,
 $\Phi_{q_1^{\beta_1}... q_{\ell}^{\beta_{\ell}}}(x) $ divides $P_{{\mathcal B}'}(x)$. Hence,  $\Phi_{q_1^{\beta_1}... q_{\ell}^{\beta_{\ell}}}(x^m) $
 divides $P_{{\mathcal B}}(x)$
 and so is $\Phi_{q_1^{\beta_1+\gamma_1}... q_{\ell}^{\beta_{\ell}+\gamma_{\ell}}} (x)$ by iteratively applying Lemma \ref{th6.1}. This implies $\B$ has $(T2)$.

\bigskip

Next we consider ${\mathcal A}$, by (\ref{eq6.2}),
$$
P_{\mathcal A}(x) = P_{{\mathcal A}'}(x^m)(1+x+...+x^{m-1}) =
P_{{\mathcal A}'}(x^m)\cdot{\prod}_{d|m,d>1}\Phi_{d}(x).
$$
We then have
$$
S_{{\mathcal A}} =  \{p^{\alpha+\gamma}:\ p^{\alpha}\in S_{{\mathcal A}'}, \   \gamma
 =\max\{j:q^{j}|m\}\}\ \cup \ \{r^j:  r  \mbox{ prime},  \ r^j|m\}.
 $$
 We can apply the same argument as the above to check that the product of the prime powers in  $S_{{\mathcal A}}$ divides $P_\A(x)$, and hence $\A$ also has $(T2)$.

  \vspace {0.2cm}

   For the last statement, we just need to note that $S_{\A} \subset  \{d: d|n\}$ by ${\mathcal A}\oplus{\mathcal B} = {\Bbb N}_n$. Hence, ${\mathcal S}$ constructed in Proposition \ref{th6.1} lies in $\frac{1}{n}{\mathbb Z}$.
\qquad$\Box$

%
%\begin{lem}\label{lem2.2}
%If ${\mathcal A} \oplus {\mathcal B} = N_n$, then
%$P_\A(x)$ and $P_\B(x)$ are co-prime.
%\end{lem}
%\proof Note that $P_\A(x)P_\B(x)=1+x+\cdots+x^{n-1}$, and the all
%roots of $1+x+\cdots+x^{n-1}$ are pair distinct. Then the result
%follows. \eproof
%
%\bigskip
%
%\noindent{\bf Proof of Theorem \ref{thm1}.} Let $\A=\SE_1\oplus
%n^{l_2}\SE_2\oplus\cdots \oplus n^{l_m}\SE_m$. According to
%Proposition \ref{prop1} and Theorem \ref{th2.1}, we can assume
%$m>1$. Let $n=\prod_{p\mid n}p^{\beta_p}$ be the decomposition of
%prime powers of $n$. Hence by definition we have
%$$S_{n^{l_i}\SE_i}=\{p^{\alpha+l_i\beta_p}: p^\alpha\in
%S_{\SE_i}\}$$
% for $1\le i\le m$ where $l_1=0$, and
% $S_\A=\bigcup_{i=1}^mS_{n^{l_i}\SE_i}$. Note that
% $$P_{n^{l_m}\SE_m}(x)=\prod_{s\in
% \widetilde{S}_{n^{l_m}\SE_m}}\Phi_s(x^{n^{l_m}}).$$
%According to Lemma \ref{Lem2.1} and \ref{lem2.2}, we have
%$\Phi_{s_1\cdots s_i}(x)\mid P_{n^{l_m}\SE_m}(x)$ for distinct prime
%powers $s_1, s_2, \ldots, s_j$ in $S_\A$. Hence $\A$ satisfies
%$(T2)$. Clearly $\A$ satisfies $(T1)$. Then by Theorem \ref{th2.1}
%the result follows.\eproof

\bigskip
\bigskip


\begin{thebibliography}{9999}

\bibitem [BG] {[BG]}
{\sc R. Bass and K. Gr\"{o}chenig}, {\it Random sampling of multivariate trigonometric polynomials}, SIAM J. Math. Anal., 36 (2004), 773-795.



\bibitem [Chr] {[Chr]}
{\sc O. Christensen}, {\it An Introduction to Frames and Riesz Bases}, Applied and Numerical Harmonic Analysis. Birkh\"{a}user Boston Inc., Boston, MA, 2003.


\bibitem [CM] {[CM]}
{\sc E. Coven and A. Meyerowitz}, {\it Tiling the integers with translates of one finite set}, J. Algebra 212(1999),
161-174.


\bibitem [DHSW]{[DHSW]}
{\sc D. Dutkay, D.G. Han, Q.Y. Sun and E. Weber}, {\it On the Beurling dimension of exponential frames},
 Adv. Math., 226 (2011), 285-297.

\bibitem [DHW]{[DHW]}
{\sc D. Dutkay, D.G. Han, and E. Weber}, {\it Bessel sequence of exponential on fractal measures}, J. Funct. Anal., 261 (2011), 2529-2539.



\bibitem [DJ] {[DJ]}
{\sc D. Dutkay and P. Jorgensen }, {\it Quasiperiodic spectra and orthogonality for iterated function system measures.}, Math. Z, 261 (2008), 373-398.


\bibitem [Fu]{[Fu]}
{\sc B. Fuglede}, {\it Commuting self-adjoint partial differential operators and a group theoretic problem}, J. Funct. Anal., 16 (1974), 101-121.

\bibitem [G] {[G]}
{\sc K. Gr\"{o}chenig}, {\it Foundations of time-frequency analysis}, Applied and Numerical Harmonic Analysis. Birkh\"{a}user, Boston, Basel, Berlin, 2001.

 \bibitem [H] {[H]}
{\sc C. Heil}, {\it A Basis Theory Primer, Expanded edition.}, Applied and Numerical Harmonic Analysis. Birkh\"{a}user Boston Inc., Boston, MA, 2011.



\bibitem [HL] {[HL]}
{\sc T.Y. Hu and K.S. Lau}, {\it Spectral property of the Bernoulli convolutions.}, Adv. Math., 219 (2008), 554-567.


\bibitem [JP] {[JP]}
{\sc P. Jorgensen and S. Pedersen}, {\it Dense analytic subspaces in fractal $L^2$ spaces.}, J. Anal. Math.,
 75 (1998), 185-228.



\bibitem [KM] {[KM]}
M. Kolountzakis and M. Matolcsi, {\it  Tiles with no spectra}, Focum Mathematicum, 18 (2006), 519 - 528.

\bibitem [L]{[L]}
{\sc C.K. Lai}, {\it On Fourier frame of absolutely continuous measures}, J. Funct. Anal., 261 (2011), 2877-2889.

\bibitem [{\L}a] {[La]} {\sc I. {\L}aba}, {\it The spectral set conjecture and multiplicative properties of roots
 of polynomials}, J. London  Math. Soc., 65 (2001), 661-671.

\bibitem [{\L}aW1] {[LaW1]}
{\sc I. {\L}aba and Y. Wang}, {\it On spectral Cantor measures}, J.
Funct. Anal., 193 (2002), 409 - 420.


\bibitem [{\L}aW2] {[LaW2]}
{\sc I. {\L}aba and Y. Wang}, {\it Some properties of spectral measures}, Appl. Comput. Harmon. Anal.,
 20 (2006), 149 -157.


\bibitem [Lan] {[Lan]}
{\sc H. Landau}, {\it Necessary density conditions for sampling and
interpolation of certain entrie functions}, Acta Math., 117 (1967),
37-52.

\bibitem [Le]{[Le]}
{\sc N. Lev}, {\it Riesz bases of exponentials on multiband spectra}, Proc. Amer. Math. Soc., to appear.


\bibitem [LLR] {[LLR]}
{\sc C.K. Lai, K.S. Lau and H. Rao}, {\it  Spectral structure of digit sets of self-similar tiles on ${\Bbb R}^1$}, preprint.



\bibitem [Lo]{[Lo]}
 {\sc C. Long},  {\it Addition theorems for sets of integers}, Pacific J. of Math., 23 (1967), 107-112.



\bibitem [LS]{[LS]}
{\sc Y. Lyubarskii and K. Seip}, {\it Sampling and interpolating sequences for multi-band-limited functions and exponential bases on disconnected sets}, J. Fourier. Anal. Appl., 3 (1997), 597-615.


\bibitem [LW]{[LW]}
{\sc J. Lagarias and Y. Wang}, {\it Spectral sets and
factorizations of finite abelian groups }, J. Funct. Anal., 145 (1997), 73 - 98.

\bibitem [OS] {[OS]}
{\sc J. Ortega-Cerd\`{a} and K. Seip}, {\it Fourier frames}, Ann of Math., 155 (2002),
789-806.

\bibitem[PW]{[PW]} {\sc S. Pedersen and Y. Wang}, {\it Universal spectra, universal tiling sets and the spectral
 set conjecture}, Math Scand., 88 (2001),  246-256.


 \bibitem [St] {[St]} {\sc R. Strichartz}, {\it Convergence of mock Fourier series }, J. Anal. Math., 99 (2006), 333-353.

      \bibitem [T] {[T]} {\sc T. Tao}, {\it Fuglede's conjecture is false in 5 or higher dimensions}, Math. Res. Letter, 11 (2004), 251-258.

     \bigskip



\end{thebibliography}
\end{document}